\documentclass{amsart}
\usepackage{ifthen, citesort}
\usepackage{amssymb,amsmath}
\usepackage{epic}
\usepackage{epsfig}

\usepackage{graphicx}
\usepackage{psfrag}

\DeclareMathOperator{\divergence}{div}

\DeclareMathOperator{\graph}{graph}
\def\R{\mathbb{R}}

\def\S{\mathbb{S}}

\def\theta{\vartheta}
\def\phi{\varphi}
\def\epsilon{\varepsilon}

\def\ul#1{\underline{#1}}
\def\ol#1{\overline{#1}}

\newtheorem{theorem}{Theorem}[section]
\newtheorem{lemma}[theorem]{Lemma}
\newtheorem{corollary}[theorem]{Corollary}

\theoremstyle{definition}
\newtheorem{remark}[theorem]{{Remark}}

\numberwithin{equation}{section}
\newcommand{\abs}[1]{\left\lvert#1\right\rvert}

\begin{document}

\title{Stability of Translating Solutions to Mean Curvature Flow}

\author{Julie Clutterbuck}
\address{
 Freie Universit\"at Berlin, Arnimallee 2-6, 
  14195 Berlin, Germany}
\curraddr{}
\email{\begin{minipage}[t]{6cm}
Julie.Clutterbuck\fuhome \\ Oliver.Schnuerer\fuhome \\ Felix.Schulze\fuhome \end{minipage}}
\def\fuhome{@math.fu-berlin.de}
\thanks{The authors are members of SFB 647/3B ``{Raum -- Zeit -- Materie:
  Singularity Structure, Long-time Behaviour and Dynamics of Solutions
  of Non-linear Evolution Equations}''}

\author{Oliver C. Schn\"urer}
\curraddr{}


\author{Felix Schulze}
\curraddr{}


\date{September 2005.}

\dedicatory{}


\begin{abstract}
We prove stability of rotationally symmetric translating 
solutions to mean curvature flow. For initial
data that converge spatially at infinity to such a soliton,
we obtain convergence for large times to that soliton 
without imposing any decay rates.
\end{abstract}

\maketitle

\section{Introduction}

We consider solutions $u:\R^n\times(0,\,\infty)\to\R$, 
$n\ge2$, to the graphical mean curvature flow equation
\begin{equation}\label{mcf} \tag{MCF}
\dot u=\sqrt{1+\lvert\nabla u\rvert^2}
\divergence\left(\frac{\nabla u}
{\sqrt{1+\lvert\nabla u\rvert^2}}\right).
\end{equation}
We will consider rotationally symmetric,
strictly convex, translating solutions $U$. 
These solutions arise as parabolic rescalings of 
type II singularities \cite{HuiskenSinestrari1}.
They have constant time derivatives $\dot U>0$.
By scaling 
we may assume that $\dot U\equiv1$,
so we can write $U(x,\,t)=U(x,\,0)+t$. 
Our main theorem states that these translating solutions 
are dynamically stable. 

\begin{theorem}\label{main thm}
Let $U:\R^n\times[0,\,\infty)\to\R,\ n\geq 2,$ be an entire rotationally symmetric, 
strictly convex solution to mean curvature flow, translating with speed $\dot U\equiv1$.
Let $u_0:\R^n\to\R$ be continuous such that the distance to $U(\cdot,0)$ 
tends to zero at infinity
$$\lim\limits_{|x|\to\infty}u_0(x)-U(x,\,0)=0.$$
\par
Then there exists a function $u\in C^{\infty}\left(\R^n\times(0,\,\infty)\right)\cap C^0\left(\R^n\times[0,\,\infty)\right)$
solving \eqref{mcf} for positive times with $u(\cdot,0)=u_0$.
\par
As time tends to infinity,
$$u(\cdot,\,t)-U(\cdot,\,t)\to0$$
uniformly on $\R^n$. 
\end{theorem}

Note that no decay rate is imposed on $u_0-U(\cdot,\,0)$. 
Our convergence 
result implies that $u$ converges to a translating
solution as $t\to\infty$. Moreover, it converges 
to precisely the translating solution we perturbed
initially.

\begin{figure}[htb]\label{barrier fig}
\psfrag{a}{upper barrier $W_R^+$}   
\psfrag{b}[c]{lower barrier $W_R^-$}   
\psfrag{c}[c]{translating convex solution $U$}   
\epsfig{file=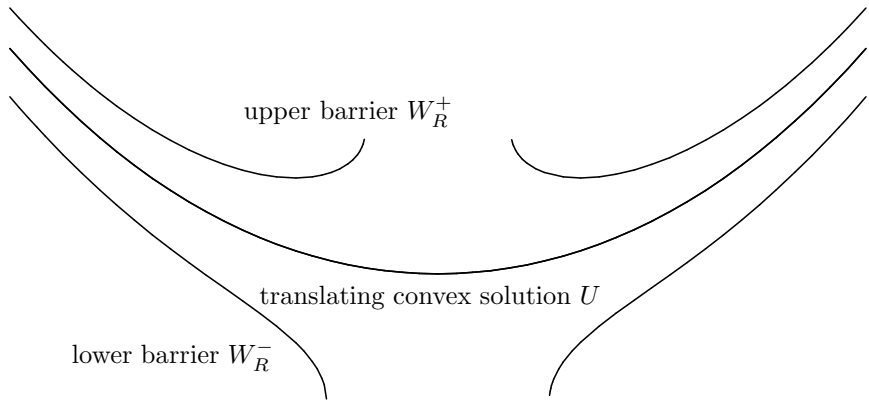, 
  width=0.9\textwidth}

\caption{Barrier Construction} 
\end{figure}

The strategy of the proof is as follows. Known results
\cite{altschuler:translating-surfaces} are easily extended
to higher dimensions and establish the existence of
a solution $U$ as in Theorem \ref{main thm}, see 
Lemma \ref{U exists}. Such rotationally symmetric translating 
solutions fulfill an ordinary differential equation. 
We will derive it in Section \ref{ode sec}. Besides these
solutions, this equation has also 
solutions which correspond to rotationally symmetric,
 graphical translating solutions which are
defined in the complement of a ball, see Figure \ref{barrier fig}
for a cross-section of both types of solutions.

We denote such a solution by $W_R^+$, if it is defined in 
$\left(\R^n\setminus{B_R(0)}\right)\times\R$ and
$\lim\limits_{|x|\downarrow R}\langle\nabla W_R^+,\,-x\rangle=\infty$.
Similarly, we call it $W_R^-$ if
$\lim\limits_{|x|\downarrow R}\langle\nabla W_R^+,\,-x\rangle=-\infty$.
Adding appropriate constants to $W_R^+$ and $W_R^-$, we see that
these solutions become asymptotic to $U$ at infinity. 

We will use these solutions as barriers. They exist in the
complement of a ball of any radius. Compare these solutions to
$n$-catenoids, $n\ge3$, which are asymptotic to hyperplanes, 
see also Appendix \ref{plane stable}.
Hyperplanes are translating solutions to mean curvature flow.
They move with velocity zero. Similarly, the translating solutions
defined in the complement of a ball are at infinity asymptotic
to entire translating solutions $U$ with $U$ as in 
Theorem \ref{main thm}.    

In order to prove stability of a translating solution,
we proceed as follows. Let $u_0$ be a graphical perturbation
of $U(\cdot,\,0)$ such that $U(\cdot,\,0)-u_0$ 
tends to zero at infinity. Then there exists a solution to
graphical mean curvature flow \eqref{mcf} which stays 
between two rotationally symmetric translating solutions. 
We obtain interior a priori estimates for this solution.
Translating solutions 
$W^+_R$ and $W^-_R$ as indicated above in 
Figure \ref{barrier fig} are then shifted vertically such 
that $W^+_R(x,t)-U(x,t)$ and $U(x,t)-W^-_R(x,t)$ tend to
$\epsilon$ for $\abs x\to\infty$. Calling these shifted solutions
again $W^+_R$ and $W^-_R$,
we may choose $R$ and $\epsilon$
such that $W^-_R(x,0)\le u(x,0)\le W^+_R(x,0)$ wherever these
functions are defined. As the gradient of $W^\pm_R(x,t)$ 
becomes unbounded for $|x|\downarrow R$, we can apply the
maximum principle and obtain that $W^-_R(x,t)\le u(x,t)\le
W^+_R(x,t)$ holds for all times. Thus $|u(x,t)-U(x,t)|$ can
attain a maximum of size larger than $2\epsilon$ only in
a bounded set. The strong maximum principle and the
a priori estimates mentioned above can be used
to show that $|u(x,t)-U(x,t)|$ will everywhere be smaller
than $2\epsilon$ for large $t$. Thus $u(x,t)-U(x,t)$
converges uniformly to zero as $t\to\infty$.

In \cite{altschuler:translating-surfaces}, Steve
Altschuler and Lani Wu have shown that entire rotationally 
symmetric translating mean curvature flow solutions exist.
They can be obtained by rescaling a type II singularity
parabolically \cite{HuiskenSinestrari1}.
Xu-Jia Wang \cite{XJWangTranslating} found other 
convex translating solutions without rotational symmetry. 
Stability of non-compact gradient K\"ahler-Ricci solitons
is investigated in \cite{OSAlbert}.
Stability of the grim reaper, a translating curve 
solving mean curvature flow, is considered in 
\cite{SmoczykStab}.
Richard Hamilton \cite{HamiltonPreprint} mentions 
non-convex, complete translating solutions to mean curvature flow. We
prove the existence of such solutions in Lemma \ref{wings exist}.

The rest of the paper is organized as follows.
In Section \ref{ode sec} we study the ordinary 
differential equation for rotationally symmetric translating
solutions. The existence of a well behaved solution is shown in
Section \ref{existence of solution}. Convergence to a translating solution 
is proved in Section \ref{conv sec}.

For the reader's convenience, we have collected
in appendices some results that we use.
We have a comparison principle in Appendix \ref{max princ sec}
and present interior estimates and an existence result 
in Appendix \ref{int ex sec}.
\par In Appendix \ref{plane stable} we show that our method 
implies directly that hyperplanes in $\R^{n+1}$, $n\ge3$, are 
stable under mean curvature flow. 
Finally, we strengthen a stability result for 
gradient K\"ahler-Ricci solitons in Appendix \ref{kr sec}.

\section{Rotationally Symmetric Translating Solutions}
\label{ode sec}
Let $V$ be a solution to graphical mean curvature flow.    For solutions that translate with speed $1$, we can write
$V(x,\,t)=V(x,\,0)+t$.  In the rotationally
symmetric case, \eqref{mcf} reduces to an
ordinary differential equation for 
$\tilde V(r)=V(x,\,0)$, where $r=\abs x$.  Writing $V(r)=\tilde V(r)$, this ordinary differential equation is 
$$1=\frac{V''}{1+V'^2}+(n-1)\frac{V'}{r},$$
where $'$ denotes derivatives with respect to $r$.
For our purposes, it will be convenient to consider the
ordinary differential equation for $\phi=V'$,
\begin{equation} \label{tr rot sym soln}
\varphi'=\left(1+\varphi^2\right)\left(1-(n-1)\frac{\varphi}r\right).
\end{equation}

Knowledge of asymptotic behavior allows us to 
find translating solutions which become close to each other at
infinity. Computer algebra calculations suggest that 
\begin{align*}
\phi=\ &\frac r{n-1}-\frac1r+(n-1)(n-4)\frac1{r^3}
-(n-1)^2\left(n^2-12n+31\right)\frac1{r^5}\\
&+(n-1)^3\left(n^3-24n^2+164n-330\right)\frac1{r^7}\\
&-(n-1)^4\left(n^4-40n^3+510n^2-2554n+4315\right)\frac1{r^9}
+O(r^{-11}).
\end{align*}
For us, the first three terms of the expansion suffice. 
\begin{lemma}\label{asymptotic expansion of phi}
For any $R>0$ and any $\varphi_0\in\R$, the boundary value problem
\begin{equation*}
\begin{cases} \phi'(r)=\left(1+\varphi^2\right)\left(1-(n-1)\dfrac{\varphi}r\right), & r\ge R,\\
\phi(R)=\phi_0, & 
\end{cases}\end{equation*}
has a unique $C^\infty$-solution $\phi$ on $[R,\infty)$. Moreover, as $r\to\infty$, we have the asymptotic expansion
\begin{equation} \varphi(r)= \frac r{n-1} -\frac1r 
+ O(r^{-2}).
\label{expansion for varphi} \end{equation}
\end{lemma}
\begin{proof}
We show that a solution exists for all $r>R$.   

Assume $1-\phi(n-1)/r$ is negative.  
Then $\phi'<0$, that is, $\phi$ is decreasing as a function of $r$.
Therefore, for sufficiently large $r$, $1-\phi(n-1)/r$ becomes positive and remains so.  
We may therefore assume that $\phi(r)<r/(n-1)$ for $r>\tilde R\ge R$.
As $\phi'>0$ in this region, $\phi$ is bounded from below.  It follows
that $\phi$ cannot become infinite for finite $r$, and we obtain existence for all $r\ge R$.  

Next, we examine the asymptotic behaviour of this solution.  
We claim that for every given $\epsilon>0$ and $r_0>R$, there exists $r_1>r_0$ such that 
$$\phi(r_1)\ge\frac{r_1}{n-1}(1-\epsilon).$$ 
If this is not the
case, we estimate $\phi'(r)\ge(1+\phi^2)\epsilon$ for $r>r_0$. 
This is impossible,
however, as our solution $\phi$ exists for $r\in[R,\infty)$.

Now observe that the function $\zeta=\frac{r}{n-1}(1-\epsilon)$ satisfies
$$ \zeta^\prime \leq \big(1+\zeta^2\big)\Big(1-\frac{n-1}{r}\zeta\Big)$$
for $r$ sufficiently large. Thus enlarging $r_1$, if necessary,
$$\frac{r}{n-1}(1-\epsilon)\leq\phi(r)\leq\frac{r}{n-1}$$ 
for $r>r_1$. Since this works for every $\epsilon>0$ we obtain
$$ \phi(r)=\frac{r}{n-1}+o(r)\ .$$

We can write $\varphi=r/(n-1)+\psi$ where $\psi$ is sublinear, so that $|\psi(r)|< cr$ for all $c>0$ and sufficiently large $r$.  This $\psi$ satisfies 
\begin{equation*}
\psi'=-\frac{(n-1)}r\psi\left(1+\left(\frac r{n-1}+\psi\right)^2\right)-\frac1{n-1}.
\end{equation*} 
Observe that $\psi$ is non-positive for $r$ sufficiently large. 
We now show that $\psi\rightarrow 0$ as $r\rightarrow\infty$.  Fix $\epsilon>0$.
Consider points $r>r_2$ where $\psi(r)\le -\epsilon$.   By the sublinearity of $\psi$, we can also assume that $-r/(2n-2)<\psi(r)$.  Then
\begin{equation*}
\psi'(r)\ge \frac{\epsilon(n-1)}{r}\left[1+\frac{ {r}^2}{4(n-1)^2}\right]-
\frac{1}{n-1}\ge c >0
\end{equation*}
for $r_2$ chosen large enough. Thus $\psi(r)\geq -\epsilon$ for $r$ large
enough.

Now set $\lambda(r)=r\psi(r)$.  We will show that $\lambda\rightarrow -1$.  Since $\psi=\lambda/r\rightarrow 0$, for all $\mu>0$ and sufficiently large $r$ we have $|\lambda(r)|\le \mu r$.   

Suppose $\lambda(r)\ge-1+\epsilon$ for some $\epsilon>0$ and
$r\geq r_3$.  Then \begin{equation*}
\lambda'(r)=-\frac r {n-1}\left[ 1 + \lambda + 2(n-1)\frac{\lambda^2}{r^2}\right]-(n-2)\frac\lambda r - (n-1)\frac{\lambda^3}{r^3}, \end{equation*}
so 
\begin{equation*}
\lambda'(r) \le -\frac{ r} {n-1}( 1 + \lambda)+ \mu(n-2)+\mu^3(n-1)
\le -c 
\end{equation*} for some $c>0$ and sufficiently large $r_3$.  We obtain that $\lambda(r)\le -1+\epsilon$ for $r$ sufficiently large.
Similarly if we assume $\lambda(r)\le -1-\epsilon$, we have 
\begin{equation*}
\lambda'(r)\ge -\frac{ r} {n-1}\left[ -\epsilon+ 2\mu^2(n-1) \right] - \mu(n-2)-\mu^3(n-1)\ge c 
\end{equation*} 
for some $c>0$ and sufficiently large $r_3$.  As above, we get $\lambda(r)\ge -1-\epsilon$ for $r$ sufficiently large.  We conclude that $\lambda\rightarrow -1$.  

Now set $\lambda(r)=-1+\eta(r)/r^2$.  We will show that $\eta\rightarrow (n-4)(n-1)$.   Since $\eta/r^2\rightarrow 0$, for all $\mu>0$ and sufficiently large $r$ we have $|\eta(r)|\le \mu r^2$. 

Suppose $\eta(r)>(n-4)(n-1)+\epsilon$ for $r>r_4$ and $\epsilon>0$.  As
\begin{align*}
\eta'(r)&=\frac{r}{n-1}\left[ (n-4)(n-1)-\eta\right]+\frac1r\left[(8-n)\eta+ (n-1)\right]
\\&\phantom{+++++}
-\frac\eta{r^3}\left[2\eta+3(n-1)\right]+3(n-1)\frac{\eta^2}{r^5}-(n-1)\frac{\eta^3}{r^7},
\end{align*}
we have 
\begin{align*}
\eta'(r)&\le \frac{r}{n-1}\left[ (n-4)(n-1)-\eta\right]+|8-n|\mu r + \frac{n-1}r 
\\&\phantom{+++++}
-\frac{\eta}{r^3}3(n-1)+3(n-1)\frac{\eta^2}{r^5}-(n-1)\frac{\eta^3}{r^7} \\
&\le \frac{r}{n-1}\left[ -\epsilon + \mu |8-n|(n-1)\right] + O(1)\\
&\le -c, 
\end{align*}
for $r_4$ large enough, if we choose $\mu$ small compared to
$\epsilon$. We obtain that $\eta(r)\le(n-4)(n-1)+\epsilon$ for $r$ sufficiently large.  
Assume now that  $\eta(r)< (n-4)(n-1)-\epsilon$ for $r>r_4$.  
Then
\begin{align*}
\eta'(r)&\ge \frac{r}{n-1}\left[ (n-4)(n-1)-\eta\right]-|8-n|\mu r + \frac{n-1}r 
\\&\phantom{+++++}
-2\mu^2r-\frac{3(n-1)\mu}{r}+3(n-1)\frac{\eta^2}{r^5}-(n-1)\frac{\eta^3}{r^7} \\
&\ge \frac{r}{n-1}\left[ \epsilon - \mu |8-n|(n-1)-2(n-1)\mu^2\right] + O(1)\\
&\ge c
\end{align*}
for small $\mu$ and large $r_4$.  We conclude that $\eta$ converges to $(n-4)(n-1)$.
\end{proof} 

\subsection{Existence of convex, rotationally symmetric translating solutions}
We mention, and slightly extend, a result of Altschuler and Wu:
\begin{lemma}\label{U exists}
There exists an entire rotationally symmetric, strictly convex graphical solution to mean curvature flow, $U:\R^n\times[0,\infty)\to\R$, $n\ge 2$, translating with speed $1$.
We have the following asymptotic expansion as $r$ approaches infinity:  
\begin{equation*} U(r,t)=t+\frac {r^2}{2(n-1)} -\ln r +O(r^{-1}).
\end{equation*}
\end{lemma}

\begin{proof}
The existence of such solutions was shown by Altschuler and Wu in \cite{altschuler:translating-surfaces} for $n=2$; their argument is in fact valid for $n\ge 2$.   They also showed that such solutions (which we denote by $U$) are asymptotic to the paraboloid $r^2/2(n-1)$.  The finer asymptotic behaviour of $U$ is a direct consequence of Lemma \ref{asymptotic expansion of phi}.
\end{proof}

\subsection{Existence of ``winglike'' translating solutions}
\label{barr exist}

\begin{figure}[htb]
\epsfig{file=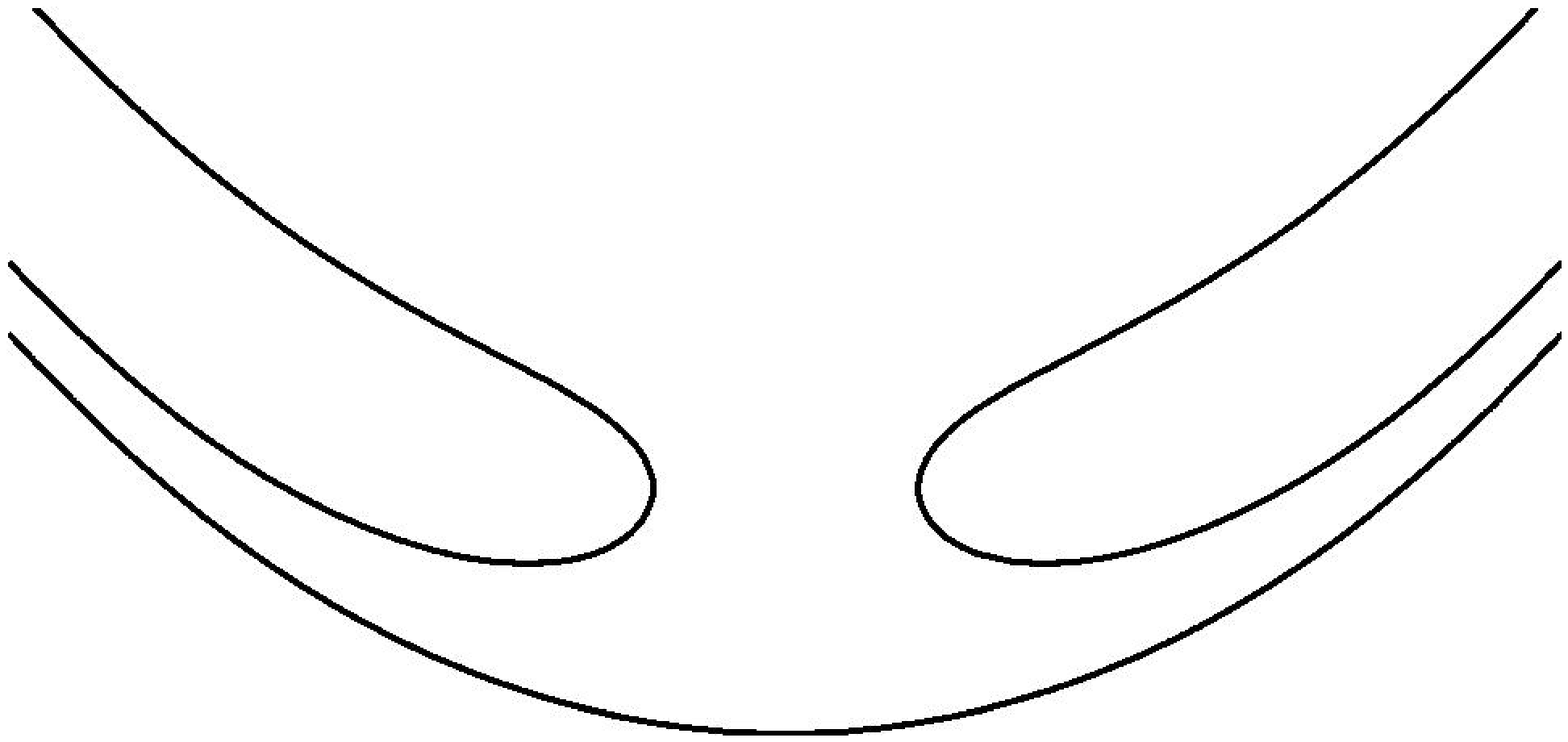, 
  width=0.9\textwidth}
\caption{Wings} 
\label{wings}
\end{figure}
We will prove that non-convex translating
solutions, made up of the union of two graphs as shown in Figure \ref{wings}, exist.

\begin{lemma}\label{wings exist}
For every $R>0$, there exist rotationally symmetric graphical solutions  to mean curvature flow, $W^+_R,W^-_R:\R^n\setminus B_R\times[0,\infty)\to\R$,  $n\ge 2$, translating with speed $1$.  
We have the following asymptotic expansion as $r$ approaches infinity:  
\begin{equation*} W^\pm_R(r,t)=t+\frac {r^2}{2(n-1)} -\ln r +O(r^{-1})+ C^\pm.
\end{equation*}
Moreover, the union of these graphs forms a complete non-convex translating solution to  mean curvature flow.
\end{lemma}
\begin{proof}
Consider a hypersurface that is invariant under rotations
around the $e_{n+1}$-axis.
Assume that this hypersurface translates with velocity
$1$ in direction $e_{n+1}$. 
At points where the tangent
plane is not orthogonal to $e_{n+1}$, we can locally
represent it as a graph $h$ over the $e_{n+1}$-axis 
$$\bigcup\limits_{x^{n+1}}h\left(x^{n+1},\,t\right)\cdot\S^{n-1}
\times\{x^{n+1}\}\subset\R^n\times\R.$$
Observe that
$h\left(x^{n+1},\,t\right)=h\left(x^{n+1}-t,\,0\right)$.  
We make the local ansatz $$x^n=\left(h^2\left(x^{n+1}-t,0\right)-\sum_{j=1}^{n-1}\left(x^j\right)^2\right)^{\frac 12}$$
and, after a tedious calculation, obtain the following 
ordinary differential equation for $h(\cdot,\,t)$
at any fixed time $t$
$$\left(\frac{n-1}h-h'\right)\left(1+h'^2\right)=h'',$$
where $'$ denotes differentiation with respect to $x^{n+1}$.
For the rest of the proof we will suppress the $t$-argument.

Fix $y_0\in\R$;
a different choice of $y_0$ corresponds to translating the
hypersurface in the $e_{n+1}$-direction. 
Starting with $h'(y_0)=0$, $h(y_0)=R$,
we obtain a strictly convex solution $h$ in a small
interval around $y_0$. Returning to our original
coordinate system and using Lemma \ref{asymptotic expansion of phi},
we can extend both branches all the way to infinity. 

We denote the upper branch by $W^-_R$ and the lower branch
by $W^+_R$. (This unintuitive naming will make sense when we
use $W^+_R$ as an upper barrier.) 
We reintroduce the time dependence according to
$$W^\pm_R(x,t):=W^\pm_R(x)+t.$$
The asymptotic behavior follows from Lemma \ref{asymptotic 
expansion of phi}. 
\end{proof}

\section{Existence of a Solution}\label{existence of solution}

\begin{theorem}\label{existence of nice solutions}
Let $u_0:\R^n\to\R$ be a continuous function with $|u_0-U(x,0)|\le d$.
Then there exists a solution $u\in C^\infty(\R^n\times(0,\infty))\cap
C^0(\R^n\times[0,\infty))$ to graphical mean curvature flow
\eqref{mcf}  with initial data $u_0$.  Furthermore 
$$\|u(\cdot,t)-U(\cdot,t)\|_{C^0} $$
is uniformly bounded for all $t\in [0,\infty)$. All derivatives of
$u-U$ are bounded locally uniformly in $x$, uniformly in $t\in[\varepsilon,\infty)$.
\end{theorem}
If $u_0$ is more regular, we can apply Theorem \ref{viscosity comparison} to obtain a uniqueness result.

\begin{proof}
We will need a local oscillation bound in order to prove uniform $C^k$ a priori
estimates. Therefore, we use Theorem \ref{thm:existDirichlet1} in such a way that we get
the oscillation bound already for the approximating solutions.

Our strategy is similar to \cite{EckerHuiskenInvent}. Our initial data $u_0$ is at finite
vertical distance from the convex rotationally symmetric translating solution $U(\cdot,0)$.
Let $\ul U(x,t):=U(x,t)-d-2$  and $\ol U(x,t):=U(x,t)+d+1$.  
Let $R>0$.
Define a continuous function $u_0^R:\ol B_R\to\R$ such that $u_0^R=u_0$ 
in $B_{R/2}$ and $u^R_0(x)=U(x,0)$ for all $x \in \partial B_R$ and
$|u_0^R(x)-U(x,0)|\leq d$ for all $x\in B_R$.
Theorem \ref{thm:existDirichlet1}
yields that there exists a solution $u^R$ to graphical mean curvature flow with 
$u^R=u_0^R$ on $B_R\times\{0\}$ and $u^R=u_0^R$ on $\partial B_R\times[0,1]$.
According to the definition of $\ol U$ and $\ul U$, we obtain for 
$(x,t)\in\partial B_R\times[0,1]$  
$$\ul U(x,t)<u^R(x,t)<\ol U(x,t).$$
We apply the standard maximum principle and deduce that
$$\ul U(x,t)<u^R(x,t)<\ol U(x,t).$$
Note that upper and lower bounds do not depend on $R$,  and so this gives us a locally uniform in $x$, uniform in $t\in[0,1]$, oscillation bound.   Applying Theorem \ref{int est} provides locally uniform in $x$, uniform in $t\in[\epsilon,1]$, $\epsilon>0$, $C^\alpha$-bounds on $u^R$:
\begin{equation*}
|D^\alpha u^R(x,t)|\le c(d,\epsilon,n,\alpha)
\end{equation*}
for $|x|<R/2$.    We let $R\to\infty$ and apply Arzel\`a-Ascoli to find a solution $u:\R^n\times(0,1]\to\R$.  

The initial data is attained for $t\downarrow0$ via the local argument used in \cite{AndrewsClutterbuck}, which survives the limiting process.  

A smooth solution for all positive times is obtained by iterating the above process.  By the interior estimates of \cite{EckerHuiskenInvent}, this solution is smooth for positive times.
We call it $u$.  It is not immediate that this solution has uniform bounds on derivatives as $t\rightarrow \infty$. 

We may now apply the comparison principle, Theorem \ref{viscosity comparison}, and deduce that \begin{equation*}
\ul U (x,t)\le u(x,t) \le \ol U(x,t) \quad\text{ for all }(x,t)\in\R^n\times[0,\infty),
\end{equation*} 
so
the oscillation of $u(x,t)-U(x,t)$  is uniformly bounded in $x$ and  $t\in[0,\infty)$.
 As above, this implies that we have bounds on $u(x,t)-U(x,t)$ and all derivatives, locally uniformly in $x$, uniformly in $t\in[\epsilon,\infty)$.  
\end{proof}

\begin{remark}
Our stability result does not require that we have a solution as
constructed in Theorem \ref{existence of nice solutions}. Suppose we 
have a $C^{2,1}(\R^n\times(0,\infty))\cap
C^0(\R^n\times[0,\infty))$-solution to graphical mean curvature flow
\eqref{mcf} with $u(x,0)=u_0(x)$ for $x\in\R^n$ and $u_0$ as in
Theorem \ref{main thm}.
Then the following lemma (and a similar estimate from below)
ensure that the given solution grows at most quadratically on bounded
time intervals.
Thus we can apply the maximum principle, Theorem \ref{viscosity
comparison} and the rest of the argument works also for
these solutions.
\end{remark}

\begin{lemma}
Let $u(x,t)$ define an entire graph moving under mean curvature flow that initially lies beneath a parabola $C|x|^2$, $C>0$.  Then $u(x,t)\le C|x|^2+ 2Cn\tau$ for $0\le t \le \tau$.  
\end{lemma} 
\begin{proof}
Fix a point $x$ with $|x|=r$, and a time $\tau$.  The sphere with radius $R$ satisfying $R^2=2n\tau + r^2+ (2C)^{-2}$ centred over the origin at height $h=(2C)^{-1}+ C(2n\tau+r^2)$ lies above our initial graph $u(\cdot,0)$.   Under mean curvature flow, the centre of the sphere remains fixed and the radius at time $t$ is $\sqrt{R^2-2nt}$.
The sphere remains above $\graph u$ for the duration of its existence, up until $T=R^2/(2n)=\tau+(2n)^{-1}(r^2+(2C)^{-2})$.  In particular it still exists at time $\tau$, at which time its height at radius $r$ is given by $Cr^2+ 2Cn\tau$.    
\end{proof}

\begin{remark}
In the above proof, we cannot directly apply the comparison principle Theorem \ref{viscosity comparison} as we do not have the growth bounds required.  On the other hand, we can compare with spheres as they are compact.
\end{remark}

\section{Convergence to a translating solution}\label{conv sec}

We first see that our wings act as barriers for $u$.
\begin{lemma}\label{wing barrier}
Let $\epsilon>0$, $R>0$.
Let $u_0:\R^n\to\R$ be a continuous function such that
$$|u_0(x)-U(x,0)|\le\epsilon\quad\text{for }|x|>R.$$
Let $u$ be a solution as in Theorem \ref{existence of nice solutions}.
Let $W^+_R$ and $W^-_R$ be ``$\epsilon$-shifted half-wing solutions'' as 
in Section \ref{barr exist}.
More precisely, we assume that $|\nabla W^+_R|$ and $|\nabla W^-_R|$
become unbounded for $|x|\downarrow R$ and
$$\lim\limits_{|x|\to\infty}W^\pm_R(x,0)-U(x,0)=\pm\epsilon.$$
Then we obtain that 
$$W^-_R(x,t)\le u(x,t)\le W^+_R(x,t)\quad\text{for }|x|>2R$$
is preserved for all times. 
\end{lemma}
\begin{proof}
We will only prove the upper bound $u(x,t)\le W^+_R(x,t)$. The lower
bound is obtained similarly. 
In order to use Corollary \ref{viscosity wt bdry}, we need to show that
\begin{equation}\label{str inequ}
u(x,t)<W^+_R(x,t)\quad\text{for }|x|=R.
\end{equation}
For small times, we may consider small spheres centred at
$(x,W^+_R(x,0))$ for $|x|=R$. These serve as barriers and show that
\eqref{str inequ} is preserved for a small time interval. 
We also obtain that
$$u(x,t)<W^+_R(x,t)\quad\text{for }|x|=R+\delta,$$
as long as $R>\delta>0$ is chosen small enough. 
This inequality holds up to some small time $t_0>0$. 

According to our a priori estimates, we get uniform 
estimates on $Du(x,t)$ for $(x,t)\in B_{2R}(0)\times[t_0,\infty)$.
Assume that $\delta$ is so small that we have
$\langle DW^+_R(x,t),x\rangle<0$ and 
$$|DW^+_R(x,t)|>2|Du(x,t)|\quad\text{for }
(x,t)\in\partial B_{R+\delta}(0)\times[t_0,\infty).$$

Our result 
follows if $u(x,t)<W^+_R(x,t)$ for $|x|=R+\delta$ and all times $t\ge 0$. If this is
not the case, there is a first time $t_1>t_0>0$ such that
$u(x_1,t_1)=W^+_R(x_1,t_1)$ for some $|x_1|=R+\delta$. 
As $|Du(x_1,t_1)|<|DW^+_R(x_1,t_1)|$ and
$\langle DW^+_R(x_1,t_1),x\rangle<0$, we deduce that
there is some $x_2$ with $|x_2|>R+\delta$ such that
$$u(x_2,t_1)>W^+_R(x_2,t_1).$$

As $t_1>0$ is minimal with that property, there exists
$0<t_2<t_1$ such that $u(x,t_2)<W^+_R(x,t_2)$ on $\{|x|=R+\delta\}$
and $u(x_2,t_2)>W^+_R(x_2,t_2)$. This contradicts 
Corollary \ref{viscosity wt bdry}.    Therefore \eqref{str inequ} is established, and with the application once again of Corollary \ref{viscosity wt bdry},  our claim follows. 
\end{proof}

Thus, we obtain that $u$ and $U$ are close to each other at infinity.  
Next, we show that $u$ becomes everywhere close to $U$ for large times.
 
\begin{lemma}\label{goes below eps}
Let $\epsilon>0$, $R>0$.
Let $u_0:\R^n\to\R$ be a continuous function such that
$$|u_0(x)-U(x,0)|\le\epsilon\quad\text{for }|x|>R.$$
Let $u$ be a solution as in Theorem \ref{existence of nice
  solutions}. Then
$$|u(x,t)-U(x,t)|\le2\epsilon$$
for $t$ sufficiently large. 
\end{lemma}
\begin{proof}
Let $W^+_R$ and $W^-_R$ be as in Lemma \ref{wing barrier}. Then 
\begin{equation}\label{wing bar eqn}
W^-_R(x,t)\le u(x,t)\le W^+_R(x,t)
\end{equation}
for $|x|>2R$.
We only show that $u(x,t)\le U(x,t)+2\epsilon$
for $t$ sufficiently large; the lower bound is obtained similarly.
Define
$$\Omega_t:=\{x\in\R^n:u(x,\,t)-U(x,\,t)>2\epsilon\}.$$
Note that $\Omega_t$ is uniformly bounded in $t$. 
We claim that $\Omega_t=\emptyset$
if $t$ is sufficiently large. As $\Omega_t$ is precompact,
$u(\cdot,\,t)-U(\cdot,\,t)$ attains its maximum somewhere
in $\Omega_t$. As in \cite[Theorem 17.1]{trudinger}, the   difference 
$w(x,\,t):=u(x,\,t)-U(x,\,t)$ fulfills a parabolic
equation of the form $\dot w=a^{ij}w_{ij}+b^iw_i$. 
If $\Omega_t\neq\emptyset$, $w(\cdot,t)$ attains its maximum, after which time
 we use the strong maximum principle to
 deduce
that this maximum is strictly decreasing in time. If 
$\Omega_t\neq\emptyset$ for all $t>0$, we define
$w_k(x,\,t):=w(x,\,t+t_k)$, 
$w_k(x,\,t):\R^n\times[-t_k,\,\infty)\to\R$, for
a sequence $t_k\to\infty$. Then a subsequence of $\lbrace w_k\rbrace$
converges locally uniformly in any $C^{l,\,l/2}$-norm
to a function $w^\infty:\R^n\times\R\to\R$ such that 
$w^\infty+U$ evolves according to \eqref{mcf}. As 
$\sup\limits_{\R^n}u(\cdot,\,t)-U(\cdot,\,t)$ is
strictly decreasing to a positive constant (again by the strict maximum principle) and 
$\Omega_t\neq\emptyset$, we deduce
that $\sup\limits_{\R^n}w^\infty(\cdot,\,t)$ is time-independent
and not smaller than $2\epsilon$. For each time $t$, 
the supremum is attained somewhere. Thus the strong maximum 
principle implies that $w$ equals a constant which is 
not smaller than $2\epsilon$, contradicting
\eqref{wing bar eqn}.
We obtain that $u(\cdot,\,t)-U(\cdot,\,t)< 2\epsilon$
for sufficiently large values of $t$. 
\end{proof}

Combining the results of this section with the existence results of Section \ref{existence of solution} gives the proof of Theorem \ref{main thm}.  

\begin{appendix}
\section{Maximum principle via viscosity solutions}
\label{max princ sec}
A special case of a comparison principle by Guy Barles, Samuel Biton, Mariane Bourgoing and Olivier Ley in \cite{BarlesetalUniqueness} is as follows:
\begin{theorem} \label{viscosity comparison}
Let $u_1,u_2:\R^n\times[0,T]\rightarrow\R$ be viscosity solutions to graphical mean curvature flow \eqref{mcf} with at most polynomial growth, that is,
\begin{equation*}
\frac{u_i(x,t)}{1+|x|^k}\to 0 \quad\text{ as }|x|\to\infty, \text{ uniformly in }t.
\end{equation*}
Let either $u_1(\cdot,0)$ or $u_2(\cdot,0)$ be locally Lipschitz continuous and fulfill 
\begin{equation*}
|Du_i(x,0)|\le C (1+|x|^\nu) \quad \text{ for almost all } x\in\R^n,
\end{equation*}
where $\nu<\left(1+\sqrt 5\right)/2$.  

If $u_1(x,0)\le u_2(x,0)$ then $u_1\le u_2$ in $\R^n\times[0,T]$.
\end{theorem}

 By direct inspection of the proof in \cite{BarlesetalUniqueness} we see that 
this result also holds if we replace $\R^n$ by $\R^n\setminus K$, $K$
a compact set, if $u_1(x,t)<u_2(x,t)$ for $x\in\partial K$.   

\begin{corollary} \label{viscosity wt bdry}
Let $u_1,u_2$ fulfill the assumptions of Theorem \ref{viscosity comparison} with 
$\R^n$ replaced by $\R^n\setminus K$ for some compact set $K$. 
If 
$$u_1<u_2$$
on $\partial K\times[0,T]$  and
$$ |u_i(x,t)|\leq C\big(1+|x|^{\nu+1}\big)$$
for all $(x,t)\in (\R^n\setminus K)\times [0,T]$, where $\nu <
\left(1+\sqrt{5}\right)/2$, we obtain that 
$$u_1\le u_2$$
on $(\R^n\setminus K)\times[0,T]$.
\end{corollary}

For positive times the solution $u$  obtained in Theorem \ref{existence of nice solutions} is a classical solution to mean curvature flow \eqref{mcf}.  It is continuous up to $t=0$.  As the locally uniform limit of a sequence of classical (and hence viscosity) solutions $\lbrace u(x,t+1/i)\rbrace $, $u$ itself is a viscosity solution for $t\in[0,T]$.

\section{Interior estimates}\label{int ex sec}
We use the following special case of a result in \cite{AndrewsClutterbuck}:

\begin{theorem}\label{thm:existDirichlet1}
 Let $u_{0}$ be a continuous function on 
$\bar B_R $ which equals a constant $b$ on $\partial B_R$.  Then there exists a unique $u\in C^{2+\alpha,1+\alpha/2}_{\text{loc}}(\bar B_R\times(0,\infty))\cap C^{0}(\bar B_R\times[0,\infty))$ satisfying graphical mean curvature flow for positive times with $u(x,0)=u_{0}(x)$ for all $x\in\bar B_R$ and $u(x,t)=b$ for $x\in\partial B_R$.  
\end{theorem}

For graphical solutions to mean curvature flow with locally in $x$ and $t$ uniformly bounded 
oscillation, there is the following interior estimate, see \cite{EckerHuiskenInvent,JCint,ColdingMinicozzi}:
\begin{theorem}\label{int est}
Let $u$ be a smooth solution to graphical mean curvature flow in
$B_R(0)\times[0,T]$ with oscillation bound $M$. Then we have
$$|Du(0,t)|\le c(R,M,t,n),$$
for $0<t<T$. We also have 
$$|D^\alpha u(0,t)|\le c(R,M,t,n,\alpha),$$
where $0<t<T$ and $\alpha$ is a multi-index denoting a combination of spatial and temporal derivatives. 
\end{theorem}

\section{Stability of the Plane}\label{plane stable}
Our method extends directly to the solution $U(x,\,t)\equiv0$, where
$(x,\,t)\in\R^n\times\R$ and $n\ge3$. An $n$-catenoid, see e.\,g.\
\cite{KaabachiPacard}, is a minimal hypersurface with two 
ends asymptotic to two parallel hyperplanes. Shifted and
rescaled appropriately, these $n$-catenoids act as barriers
and imply stability of $\R^n\subset\R^{n+1}$, $n\ge3$, under
mean curvature flow.
\begin{theorem}
Let $u_0:\R^n\to\R$ be a continuous function decaying at 
infinity,
$$\lim\limits_{\abs x\to\infty}u_0(x)=0.$$ 
Then there exists a solution
$u\in C^\infty\left(\R^n\times(0,\infty)\right)\cap C^0\left(\R^n\times[0,\infty)\right)$
to \eqref{mcf} for positive times with $u(\cdot,0)=u_0$.
For $t\to\infty$, $u(x,\,t)$ converges uniformly to zero.
\end{theorem}
Amusingly, we use $n$-catenoids, which are unstable
minimal hypersurfaces \cite{CaoShenZhu},
in order to prove stability.


\section{Stability of Gradient K\"ahler-Ricci Solitons}
\label{kr sec}

In \cite{OSAlbert}, Albert Chau and the second author 
proved stability for gradient K\"ahler-Ricci flow solitons.
These solitons are analogous to the translating solutions
considered here. 

The main theorem of that paper, in which a decay rate is 
imposed, may be extended. We can drop the decay condition
(2) in \cite[Theorem 1.2]{OSAlbert} in favour of the requirement that the solution 
$u$ to K\"ahler-Ricci flow  initially  tends to zero at infinity. 

It suffices to add a small positive constant $\epsilon$ to an upper barrier 
used in the proof there and to argue as in the proof of 
Theorem \ref{main thm} above. In this way, we show that $u$ is
eventually smaller than $2\epsilon$. 
\end{appendix}

\bibliographystyle{amsplain}

\begin{thebibliography}{10}

\bibitem{altschuler:translating-surfaces}
Steven~J. Altschuler and Lang~F. Wu, \emph{Translating surfaces of the
  non-parametric mean curvature flow with prescribed contact angle}, Calc. Var.
  Partial Differential Equations \textbf{2} (1994), no.~1, 101--111.

\bibitem{AndrewsClutterbuck}
Ben Andrews and Julie Clutterbuck, \emph{Time-interior gradient estimates for
  graphical an\-iso\-tropic mean curvature flows}, 2005.

\bibitem{BarlesetalUniqueness}
Guy Barles, Samuel Biton, Mariane Bourgoing, and Olivier Ley, \emph{Uniqueness
  results for quasilinear parabolic equations through viscosity solutions'
  methods}, Calc. Var. Partial Differential Equations \textbf{18} (2003),
  no.~2, 159--179.

\bibitem{CaoShenZhu}
Huai-Dong Cao, Ying Shen, and Shunhui Zhu, \emph{The structure of stable
  minimal hypersurfaces in {${\bf R}\sp {n+1}$}}, Math. Res. Lett. \textbf{4}
  (1997), no.~5, 637--644.

\bibitem{OSAlbert}
Albert Chau and Oliver~C. Schn\"urer, \emph{{S}tability of gradient
  {K}\"ahler-{R}icci solitons}, 2003, MPI-MIS Preprint 39/2003, {\tt
  http://www.mis.mpg.de/}, to appear in Comm. Anal. Geom..

\bibitem{JCint}
Julie Clutterbuck, \emph{Parabolic equations with continuous initial data},
  Ph.D. thesis, Australian National University, 2004.

\bibitem{ColdingMinicozzi}
Tobias~H. Colding and William~P. Minicozzi~II, \emph{{Sharp estimates for mean
  curvature flow of graphs.}}, J. Reine Angew. Math. \textbf{574} (2004),
  187--195.

\bibitem{EckerHuiskenInvent}
Klaus Ecker and Gerhard Huisken, \emph{Interior estimates for hypersurfaces
  moving by mean curvature}, Invent. Math. \textbf{105} (1991), no.~3,
  547--569.

\bibitem{trudinger}
David Gilbarg and Neil~S. Trudinger, \emph{Elliptic partial differential
  equations of second order}, Classics in Mathematics, Springer-Verlag, Berlin,
  2001, Reprint of the 1998 edition.

\bibitem{HamiltonPreprint}
Richard~S. Hamilton, \emph{Eternal solutions to the mean curvature flow},
  unpublished.

\bibitem{HuiskenSinestrari1}
Gerhard Huisken and Carlo Sinestrari, \emph{Mean curvature flow singularities
  for mean convex surfaces}, Calc. Var. Partial Differential Equations
  \textbf{8} (1999), no.~1, 1--14.

\bibitem{KaabachiPacard}
Saida Kaabachi and Frank Pacard, \emph{{Riemann minimal surfaces in higher
  dimensions}}, \rule{0mm}{1mm}\\ {\tt
  http://perso-math.univ-mlv.fr/users/pacard.frank/}.

\bibitem{SmoczykStab}
Knut Smoczyk, \emph{A relation between mean curvature flow solitons and minimal
  submanifolds}, Math. Nachr. \textbf{229} (2001), 175--186.

\bibitem{XJWangTranslating}
Xu-Jia Wang, \emph{{Convex solutions to the mean curvature flow}}, {\tt
  arXiv:math.DG/0404326}.

\end{thebibliography}


\end{document}